\documentclass[a4paper,12pt]{amsart}

\usepackage{amsmath,amsfonts,amssymb,amsthm}
\usepackage{url}

\newcommand{\R}{\mathbb{R}}

\usepackage{graphicx}
\graphicspath{{./pictures/}}

\begin{document}

\title[Weak instability]{Weak instability of Hamiltonian equilibria}
\author[Gaetano Zampieri]{Gaetano Zampieri}
\address{Universit\`a di Verona,
Dipartimento di Informatica\\
strada Le Grazie 15, 37134 Verona, Italy}
\email{gaetano.zampieri@univr.it}

\maketitle

\begin{abstract} This is an expository paper on Lyapunov stability of equilibria of autonomous Hamiltonian systems.
Our aim is to clarify the concept 
of \emph{weak instability}, namely instability without non-constant motions which have the equilibrium as limit point as time goes to minus infinity.  This is done by means of some examples. In particular, we show that
a weakly unstable equilibrium point can be  stable for the linearized vector field.
\end{abstract}

\section{Introduction} Stability of the equilibrium is a mathematical field more then two centuries old.
Indeed,  Lagrange   stated the celebrated Lagrange-Dirich\-let theorem in the eighteenth century, and some so called converses of that statement are still proved nowadays. 
So many mathematicians have been interested in stability   that we   refrain from
mentioning them with the exception of the most important,   Lyapunov, who defended his   doctoral thesis ``The general problem of the stability of motion'' in 1892.  The applications are also  countless 
in  mechanics and in most sciences. To start with the rich literature   on this matter, see   Arnold et al.~\cite{arnold},  Meyer et al.~\cite{meyer}, and Rouche et al.~\cite{rouche}.

Important mathematical objects related to the instability of the equilibrium are  \emph{asymptotic motions}. Before
their formal  definition, let us mention that the upper position of a simple pendulum, and zero velocity, constitute an unstable
equilibrium and its asymptotic motions are neither rotations (when the pendulum swings around and around) nor librations (when it swings back and forth), and they stay between the two behaviors.

Let us consider a smooth vector field $f$ on an open set $A\subseteq\R^N$ with an equilibrium point $\hat x\in A$, so $f(\hat x)=0$. We say
that $\phi:(-\infty,b)\to A$ is an \emph{asymptotic motion} in the past to the equilibrium point $\hat x$, if
$\phi(t)$ is a non-constant solution to the o.d.e. $\dot x=f(x)$ such that $\phi(t)\to \hat x$ as $t\to-\infty$.
 In the sequel we briefly write `asymptotic motion' instead of `asymptotic motion in the past' since we are only concerned with this
kind of asymptotic motions. Of course the existence of an asymptotic motion implies the Lyapunov instability of the equilibrium point.
The basic sufficient condition for the existence of an asymptotic motion is the presence of an eigenvalue of $f'(\hat x)$ with
strictly positive real part, see for instance Hartman~\cite{hartman} remark to Corollary 6.1, p.~243.

In this paper we focus on \emph{autonomous Hamiltonian systems} so in the sequel $N=2n$, $x=(q,p)$, $q,p\in\R^n$, and the vector field is
\begin{equation}
\left(\partial_p H(q,p), -\partial_q H(q,p)\right)
\end{equation} for some smooth  $H$ called the Hamiltonian function. 
Our aim is to clarify the concept 
of \emph{Lyapunov instability without asymptotic motions} that we briefly call \emph{weak instability}. This is done by means of some examples.

Section \ref{linear} deals with linear systems.
Of course there is a trivial situation where weak instability appears: the  free particle. 
The equilibrium is non-isolated and the eigenvalues vanish, the example can be done in one degree of freedom so in dimension 2. A more subtle instability of the equilibrium for a linear system is obtained when the eigenvalues are purely imaginary and some Jordan blocks have dimension greater than one,
of course this can happen only in dimension at least~4. The example we are going to see comes from the planar restricted 3-body problem at one of the relative equilibria, the Lagrange equilateral   points, also called  the Trojan points, at the critical Routh value of the mass ratio of the primaries.

In Section~\ref{nonlinear} we move on nonlinear systems. Their equilibria can be unstable even if we have stability
for the linearized system
as the Cherry Hamiltonian  in dimension~4 shows by means of an asymptotic motion. Cherry's system  is the third example of this paper, it was published in 1925 and, in the last 20 years, it became important in plasma physics, see Pfirsch~\cite{pfirsch} and the references therein. 

Our fourth  example, 
also in dimension~4, comes from~\cite{zampieri3} and shows that \emph{we can have weak instability
of an Hamiltonian equilibrium which is linearly stable}. 
Some  systems, 
produced by Barone-Netto and myself~\cite{zampieri1} and \cite{zampieri2}, preceded~\cite{zampieri3},  they give non-Hamiltonian examples of weak instability for linearly stable equilibria. 

Hopefully,  the concept of weak instability will 
 stimulate further researches  in stability within mathematical physics, together
 with other fresh notions like the ``weak asymptotic stability'' introduced by Ortega, Planas-Bielsa and Ratiu, see~\cite{ortega} and the references therein.

%%%%%%%%%%%%%%%%%%%%%%%%%%%%%%%
\section{Weak instability for linear  systems}\label{linear} 
\subsection{Free particle}
 Our first example is a particle on a straight line under no forces
\begin{equation}\label{HforFreeParticle}
 H(q,p)=\frac{p^2}{2},\qquad q,p\in\R.
 \end{equation}
The Hamiltonian vector field  is 
\begin{equation}\label{Hfree}
\bigl(\partial_{p} H(q,p),-\partial_{q} H(q,p)\bigr)=
(p,0).
\end{equation}
It is a linear  field with the double eigenvalue $0$.
The integral curves are
\begin{equation}
q(t)=q(0)+p(0)t\,,\qquad
p(t)=p(0).
\end{equation}
Each $(q_0,0)\in\R^2$ is an equilibrium point and its instability  can be shown by means of the sequence $(q(0),p(0))=(q_0,1/m)\to (q_0,0)$ as $m\to +\infty$. There are no asymptotic motions.

\subsection{Linearization at $L_4$}
Our second example is the quadratic part of the
Hamiltonian function of the planar restricted 3-body problem at one of the relative equilibria,
the Lagrange libration point $L_4$ at the critical Routh value of the mass ratio of the primaries. In the sequel
  $q=(q_1,q_2), p=(p_1,p_2), (q,p)=(q_1,q_2,p_1,p_2)$, and 
\begin{equation}\label{HforLagrangeL4}
 H(q,p)=\frac{1}{\sqrt 2}\det(p,q)+\frac{1}{2}|q|^2=\frac{1}{\sqrt 2}\bigl(p_1 q_2-p_2q_1\bigr)+\frac{1}{2}\bigl(q_1^2+q_2^2\bigr).\end{equation}
The Hamiltonian vector field  is
\begin{equation}\label{HvfforLagrangeL4}\begin{split}
&\left(\partial_{p_1} H(q,p),
\partial_{p_2} H(q,p),-\partial_{q_1} H(q,p),
-\partial_{q_2} H(q,p)
\right)=\\
&=\left(q_2/\sqrt{2},- q_1/\sqrt{2},-q_1+p_2/\sqrt{2},-q_2-p_1/\sqrt{2}\right),\end{split}
\end{equation}
see $H_0$ and the o.d.e. at the end of p.~256, with $\xi=q$, $\eta=p$, $\omega=1/\sqrt{2}$, $\delta=1$, and also $H_0$   at p.~258 in Meyer et al.~\cite{meyer}.

It is a linear vector field with the double eigenvalues   $\lambda=\pm i/\sqrt{2}$ and Jordan blocks
$\begin{pmatrix} \lambda &1\\0 &\lambda\end{pmatrix}$. The origin is now the unique equilibrium point.
  
The function $|q|^2$ is a first integral. Suppose the integral curve  $(q(t),p(t))\to 0$ as $t\to -\infty$, then $|q(t)|^2\equiv 0$ and  this fact further implies that $|p(t)|^2\equiv 0$, indeed for $q(t)\equiv 0$ we have
\begin{equation}
\frac{d}{dt}|p(t)|^2=2p(t)\cdot \left(-q_1(t)+p_2(t)/\sqrt{2},-q_2(t)-p_1(t)/\sqrt{2}\right)= 0.
\end{equation}
 So the integral curve is constant and we do not have
asymptotic motions to the equilibrium point.

The origin is an unstable equilibrium point as we can see with 
\begin{equation}\label{solutionforLagrangeL4}\begin{split}
q_1(t)&=\frac{1}{m}\cos\frac{t}{\sqrt{2}}\,,\qquad \quad  q_2(t)=-\frac{1}{m}\sin\frac{t}{\sqrt{2}}\,,\\
p_1(t)&=-\frac{t}{m}\cos\frac{t}{\sqrt{2}}\,,\qquad\ 
p_2(t)=\frac{t}{m}\sin\frac{t}{\sqrt{2}}\,,
\end{split}
\end{equation}
for w
$(q_1(0),q_2(0),p_1(0),p_2(0))=(1/m,0,0,0)\to 0$ as  $m\to+\infty$.

Incidentally, in connection with the nonlinear 3-body problem which has the Hamiltonian vector field defined by~\eqref{HvfforLagrangeL4} as linearization at $L_4$, the book \cite{meyer} at the end of Sec.~13.6 says that in 1977 two papers claimed
to have proved the stability of the equilibrium, however one proof is wrong and the other is unconvincing. The last sentence is: ``It would be interesting to give a correct proof of stability in this case, because the linearized system is not simple, and so the linearized equations are unstable''.

%%%%%%%%%%%%%%%%%%%%%%%%%%%%%%%
\section{Instability for linearly stable  equilibria}\label{nonlinear} 
\subsection{Cherry Hamiltonian}
Next,  the famous Cherry Hamiltonian system  shows that the equilibrium can be unstable
even if it is stable for the linearized system, briefly even if it is \emph{linearly stable}.
In  Cherry~\cite{cherry} p.~199, or in Whittaker \cite{whittaker} p.~412, we can see the Hamiltonian function $H:\R^4\to\R$
\begin{equation}\label{CherryH}
H(q,p)=\frac{1}{2}\bigl(q_1^2+p_1^2\bigr)-\bigl(q_2^2+p_2^2\bigr)+\sigma\Bigl(q_2\bigl(q_1^2-p_1^2\bigr)-2 q_1p_1p_2\Bigr).
\end{equation}
The Hamiltonian vector field, written as a column vector, is
\begin{equation}\label{HvfforCherry}
\begin{pmatrix}p_1-2\sigma q_2 p_1-2\sigma q_1p_2\\-2 p_2-2\sigma q_1p_1\\-q_1-2\sigma q_2q_1+2\sigma p_1p_2\\
2q_2+\sigma p_1^2-\sigma q_1^2  
\end{pmatrix}.
\end{equation}
The  linearized vector field  $(p_1,-2p_2,-q_1,2q_2)$ is obtained  for
 $\sigma=0$. The origin is stable for the linearized systems which consists of two  harmonic oscillators: $\ddot q_1=-q_1$, 
$\ddot q_2=-4q_2$. The eigenvalues are distinct $\pm i$, $\pm 2i$. However, the origin is Lyapunov unstable for the vector field~\eqref{HvfforCherry} whenever $\sigma\ne 0$ since it has the following
asymptotic motion  defined for $t<0$
\begin{equation}\label{as}
\begin{split}
q_1(t)&= \frac{\sin t}{\sqrt{2}\,\sigma\, t}\,,\qquad 
q_2(t)=\frac{\sin(2t)}{2\sigma\,t}\,,\\ 
p_1(t)&=\frac{\cos t}{\sqrt{2}\,\sigma\,t}\,,\qquad  
p_2(t)=-\frac{\cos(2t)}{2\sigma\, t}\,.
\end{split}
\end{equation}
\begin{figure}
  \begin{center}
\includegraphics[scale=.55]{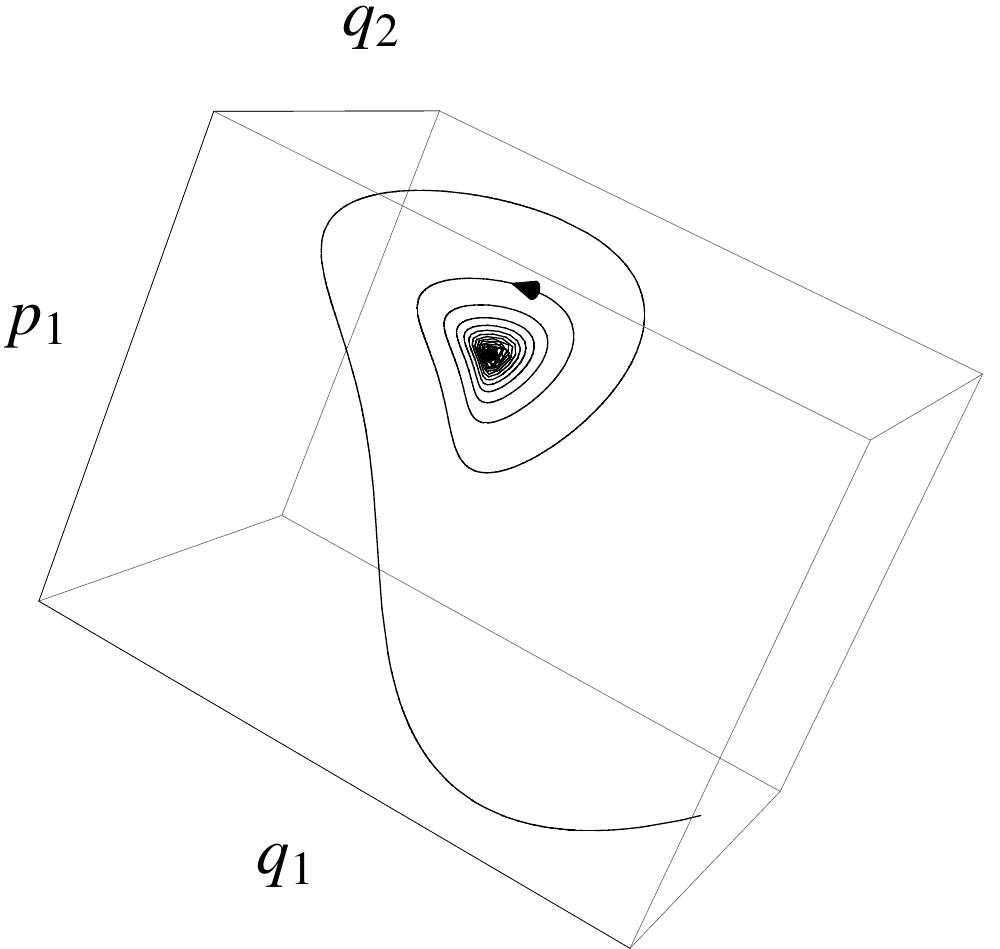}
\end{center}
\caption{Asymptotic motion for Cherry Hamiltonian}
\label{asintotico}
\end{figure}

%%%%%%%%%%%%%
\subsection{Variation-like Hamiltonian} 
Our final  example shows  that the origin is an  unstable equilibrium point which is linearly stable and
has no asymptotic motions for the system defined by
\begin{equation}\label{polynomial}
 H(q,p)=p_1p_2+q_1q_2+\sigma\, q_1^2 q_2,\qquad \sigma\ne 0.
\end{equation}
It is a particular case of the following Hamiltonian function introduced in~\cite{zampieri3}
\begin{equation}\label{HZamp}
  H(q,p)=p_1p_2+g(q_1)q_2,\quad\quad g(0)=0,\ g'(0)>0,
  \end{equation}
  where $ g\in C^1$  on a neighborhood of $0$. The Hamiltonian vector field is
\begin{equation}\label{HvfforZamp}
\begin{pmatrix}p_2\\p_1\\-g'(q_1)q_2\\-g(q_1)  
\end{pmatrix}=\begin{pmatrix}p_2\\p_1\\-g'(0)q_2\\-g'(0)q_1  
\end{pmatrix}+o\bigl(|(q,p)|\bigr).
\end{equation}
The origin is stable for the linearized system which  consists of two  harmonic oscillators: $\ddot q_1=-g'(0)q_1$, 
$\ddot q_2=-g'(0)q_2$. In this case the eigenvalues are double $\pm i\sqrt{g'(0)}$ however the Jordan blocks are one-dimensional.  

The subsystem of the first and last canonical equations 
\begin{equation}\label{subsystem}
\dot q_1=p_2,\qquad \dot p_2=-g(q_1),
\end{equation}
 separates. If we take a 
solution $(q_1(t),p_2(t))$ of this subsystem and plug $q_1(t)$  into the second and third canonical equations, we then get the  
equations of variation of \eqref{subsystem} along the solution $(q_1(t),p_2(t))$. This is why
the function in formula \eqref{HZamp} is called \emph{variation-like Hamiltonian} in the title of this subsection.

There are no asymptotic motions, indeed if the solution 
\begin{equation}
(q_1(t),q_2(t),p_1(t),p_2(t)\bigr)\to 0\qquad \hbox{\rm as}\quad t\to -\infty
\end{equation} then 
$(q_1(t),p_2(t))\equiv 0$, since the origin is a local center for \eqref{subsystem},   and this implies $(q_2(t),p_1(t))\equiv 0$ too.

\begin{figure}
  \begin{center}
\includegraphics[scale=.5]{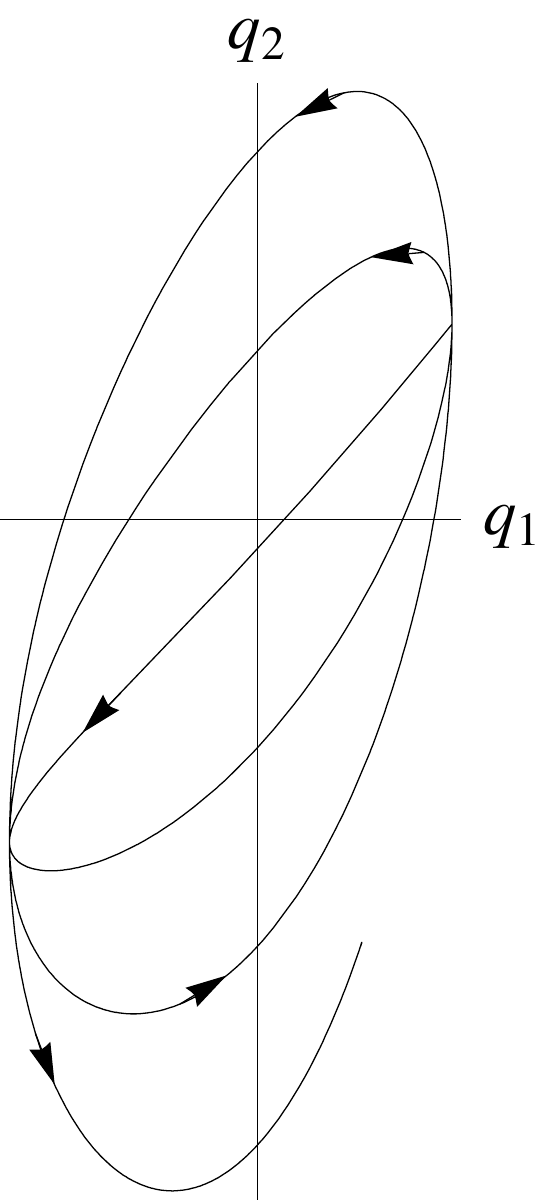}
\end{center}
\caption{Projection of an unbounded orbit for $H=p_1p_2+q_1q_2+q_1^2 q_2$}
\label{unbounded}
\end{figure}

In spite of this fact, the origin is unstable for \eqref{HvfforZamp} for most functions $g$ as above. Theorem~3.3 in
\cite{zampieri3} proves that stability is equivalent to the isochrony of the periodic solutions of the subsystem \eqref{subsystem} in a neighborhood of $0\in\R^2$, and this  implies the isochronous periodicity of all integral curves of
\eqref{HvfforZamp} in a neighborhood of $0\in\R^4$. Moreover, 
 Corollary~2.3 in \cite{zampieri3}  for a smooth $g$ provides
\begin{equation}
g'''(0)=\frac{5 g''(0)^2}{3 g'(0)} 
\end{equation}
as the simplest necessary condition for (local isochrony and then) stability.
So the choice $g(q_1)= q_1+\sigma\, q_1^2$ of the Hamiltonian~\eqref{polynomial} gives instability for all $\sigma\ne 0$. 
In Figure~\ref{unbounded} we can see the projection on the $q_1,q_2$-plane of the integral curve of the Hamiltonian vector field given by~\eqref{polynomial}. 

Finally, let us remark that the Hamiltonian \eqref{polynomial}, composed with the symplectic transformation $(Q,P)\mapsto(Q_1+Q_2,Q_1-Q_2,P_1+P_2,P_1-P_2)/\sqrt{2}$, becomes  
\begin{equation}
\frac{1}{2}\bigl(Q_1^2+P_1^2\bigr)-\frac{1}{2}\bigl(Q_2^2+P_2^2\bigr)+\frac{\sigma}{2\sqrt{2}}(Q_1+Q_2)\bigl(Q_1^2-Q_2^2\bigr)
\end{equation}
 a function with some  features in common with Cherry's Hamiltonian~\eqref{CherryH}.

%%%%%%%%%%%%%%%%%%%%%%%%%%%%%%%
\paragraph{Acknowledgments.}   The pictures  were  made using  {\it
Mathematica}  by    Wolfram Research Inc. by means of the   package
CurvesGraphics6  by
Gian\-luca Gorni  available at:
\url{http://sole.dimi.uniud.it/~gianluca.gorni/}


\begin{thebibliography}{9}


\bibitem{arnold} Vladimir I. Arnold, Valery V. Kozlov, and Anatoly I. Neishtadt, \textit{Mathematical Aspects of Classical and Celestial Mechanics}. 
Encyclopaedia of Mathematical Sciences~3. Dynamical Systems III. Springer-Verlag,
1988.

\bibitem{cherry} Thomas M. Cherry, Some examples of trajectories defined by differential equations of a generalized dynamical type, \emph{Trans. Cambridge Phil. Soc.} \textbf{23} (1925) 169--200.

\bibitem{hartman} Philip Hartmann, \textit{Ordinary Differential Equations}. Second ed.  Birkh\"auser, 1982.

\bibitem{meyer} Kenneth R. Meyer,  Glen R. Hall, and  Dan Offin, \textit{Introduction to Hamiltonian Dynamical Systems and the N-Body
Problem}. Second ed. Applied Mathematical Sciences, 90. Springer-Verlag, 2009.

\bibitem{ortega} Juan-Pablo Ortega, V\'ictor Planas-Bielsa, Tudor S. Ratiu, 
Asymptotic and Lyapunov stability of constrained and Poisson equilibria, 
\emph{J. Diff. Eqns.} \textbf{214} (2005) 92--127.

\bibitem{pfirsch} Dieter Pfirsch, Nonlinear instabilities relating the negative-energy modes, \emph{Physical Review E} \textbf{48} (1993) 1428--1435.

\bibitem{rouche} Nicolas Rouche, Patrick Habets, and Michel Laloy, 
\textit{Stability theory by Liapunov's direct method}. 
Applied Mathematical Sciences,  22. Springer-Verlag,  1977.

\bibitem{whittaker} Edmund T. Whittaker, \textit{A Treatise on the Analytical Dynamics of Particles and Rigid Bodies}. Fourth ed. Dover Publications, 1944.

\bibitem{zampieri1} Gaetano Zampieri,  and Angelo Barone-Netto,  
Attractive central forces may yield Liapunov instability, in \textit{Dynamical systems and partial differential equations}, Proceedings of the VII Elam, Editorial Equinoccio, Caracas, 1986, 
  105--112.

\bibitem{zampieri2} Gaetano Zampieri, Liapunov stability for some central
forces, \emph{J. Diff. Eqns.} \textbf{74} (1988) 254--265.

\bibitem{zampieri3} Gaetano Zampieri, Completely integrable
Hamiltonian systems with weak Lyapunov instability
or isochrony,  \emph{Commun. Math. Phys.} \textbf{303} (2011) 73--87.


\end{thebibliography}
\end{document}